\documentclass[12pt]{article}
\usepackage[reqno]{amsmath}
\usepackage{amssymb, theorem, enumerate}
\usepackage{graphicx}
\usepackage{tikz}
\usetikzlibrary{decorations.markings}
\usepackage{array}

\expandafter\ifx\csname pdfpagebox\endcsname\relax\else
\fi

\theoremstyle{change}
{\theorembodyfont{\slshape}
\newtheorem{theorem}{Theorem.}[section]

}
\theorembodyfont{\rmfamily}

\newcommand\cref[1]{Corollary~\ref{cor:#1}}

\def\proof{\noindent{{\sl Proof. }}}
\def\sqr#1#2{{\vbox{\hrule height.#2pt
    \hbox{\vrule width.#2pt height#1pt \kern#1pt
        \vrule width.#2pt}\hrule height.#2pt}}}
\def\eqed{\sqr53}
\def\qed{%
    \ifmmode\eqno\eqed
    \else\nobreak\ \hfill\eqed\medbreak\fi}

\newcommand\Zy{{\mathbf y}}

\newcommand\Zz{{\mathbf z}}
\newcommand\Zx{{\mathbf x}}

\newcommand\cx{{\mathbb C}}% complexes

\newcommand\re{{\mathbb R}}%reals

\DeclareMathOperator{\tr}{tr}
\DeclareMathOperator{\col}{col}
\DeclareMathOperator{\spec}{sp}

\newcommand\ins{D_{h}}
\newcommand\outs{D_{t}}
\newcommand\splus{S^+(U)}

\title{Quantum Walks on Regular Graphs and Eigenvalues}
\author{Chris Godsil and Krystal Guo}

\begin{document}
\maketitle

\begin{abstract}
	We study the transition matrix of a quantum walk on strongly regular graphs. It is proposed by Emms, Hancock, Severini and Wilson in 2006, that the spectrum of $S^+(U^3)$, a matrix based on the amplitudes of walks in the quantum walk, distinguishes strongly regular graphs. We find the eigenvalues of $S^+(U)$ and $S^+(U^2)$ for regular graphs and show that $S^+(U^2) = S^+(U)^2 + I$. 
\end{abstract}

\section{Introduction}

A discrete-time quantum walk is a quantum process on a graph whose state vector is governed by a matrix, called the transition matrix. In \cite{ESWH, EHSW06} Emms, Severini, Wilson and Hancock propose that the quantum walk transition matrix can be used to distinguish between non-isomorphic graphs. Let $U(G)$ and $U(H)$ be the transition matrices of quantum walks on $G$ and $H$ respectively. Given a matrix $M$, the \textsl{positive support} of $M$, denoted $S^+(M)$, is the matrix obtained from $M$ as follows:
\[ (S^+(M))_{i,j} = \begin{cases} 1 & \text{if } M_{i,j} >0\\
0 & \text{otherwise.}\end{cases}
\]

\begin{theorem}\label{su1} If $G$ and $H$ are isomorphic  regular graphs, then $S^+(U(G)^3)$ and $S^+(U(H)^3)$ are cospectral. \end{theorem}

The authors of \cite{EHSW06, ESWH} propose that the converse of Theorem \ref{su1} is also true; they conjecture that the spectrum of the matrix $S^+(U^3)$ distinguishes strongly regular graphs. After experiments on a large set of graphs, no strongly regular graph is known to have a cospectral mate with respect to this invariant. If the conjecture is true, it would yield a classical polynomial-time algorithm for the Graph Isomorphism Problem for strongly regular graphs (but there do not seem to be strong grounds for believing the conjecture). 

In this paper we will find the spectra of two matrices related to proposed graph invariant, for regular graphs. In \cite{EHSW06}, Emms et al.~compute some eigenvalues of $\splus$ and $S^+(U^2)$ but do not determine them all; for both matrices, they find the set of eigenvalues which are derived from the eigenvalues of the adjacency matrix, but do not find the remaining eigenvalues. The spectrum of $S^+(U)$ is also given in \cite{raewh11}. 

Here we will use an approach which exploits the linear algebraic properties of $\splus$ to yield a proof that the spectrum of $\splus$ is determined by the spectrum of the graph with respect to the adjacency matrix. We also completely determine the spectrum of $S^+(U^2)$ by expressing $S^+(U^2)$ in terms of $\splus$ and the identity matrix. 

\section{Preliminary Definitions}
\label{sec:Defs}

A \textsl{discrete-time quantum walk} is a process on a graph $G$ governed by a unitary matrix, $U$, which is called the \textsl{transition matrix.} For $uv$ and $wx$ arcs in the digraph of $G$, the transition matrix is defined to be:
\[
 U_{wx,uv} = \begin{cases} \frac{2}{d(v)} &\text{ if } v=w \text{ and } u \neq x ,\\
\frac{2}{d(v)} -1  & \text{ if } v=w \text{ and } u = x, \\
0 &\text{ otherwise.} \end{cases}
\]

Let $A$ the adjacency matrix of $G$. Let $D$ be the digraph of $G$ and consider the following incidence matrices of $D$, both with rows indexed by the vertices of $D$ and columns indexed by the arcs of $D$:
\[ (\ins)_{i,j} = \begin{cases} 1 &\text{if } i \text{ is the head of arc }j \\ 
0 &\text{otherwise}\end{cases}
\]
and 
\[ (\outs)_{i,j}  = \begin{cases} 1 &\text{if } i \text{ is the tail of arc }j \\ 
0 &\text{otherwise.}\end{cases}
\]
To describe the quantum walk, we need one more matrix: let $P$ be a permutation matrix with row and columns indexed by the arcs of $D$ such that,
\[P_{wx,uv}  = \begin{cases} 1 &\text{if } x=u \text{ is the tail of arc } w=v \\ 
0 &\text{otherwise.}\end{cases}
\]
Then, we see that $\ins\outs^{T} = A(G)$ and
\[
 (\outs^T\ins)_{wx,uv} = \begin{cases} 1 &\text{if } v=w,\\
0 &\text{otherwise.} \end{cases}
\]
If $G$ is regular with valency $k$, we have that \[U = \frac{2}{k} \outs^T\ins -  P.\] 

\section{Eigenvalues of $S^+(U)$}
\label{sec:su}

In this section, we will find the eigenvalues of $S^+(U)$ for a regular graph $G$. If $G$ is regular with valency $1$, then $G$ must be a matching and the spectrum of $S^+(U(G))$ is easily determined. We may direct our attention to regular graphs with valency $k \geq 2$. If $G$ is a regular graph with valency $k$ on $n$ vertices, then \[U =  \frac{2}{k} \outs^T\ins - P.\] The only negative entries have values $\frac{2}{k} -1$, for $k \geq 2$, so $S^+(U) = \outs^T\ins - P$. 

From Section \ref{sec:Defs}, we see that $\outs \outs^T = k I$ and $\ins \ins^T = k I$.  From the definition of $P$, we get that 
\[ P \ins^T = \outs^T \; \text{ and } \; P\outs^T = \outs^T
\]
Let $Q = \frac{2}{k} \ins^T \ins -I$. Then, $ Q^2 = I$ and we can write $\splus$ as:
\[ S^+(U) = \outs^T\ins - P = P ( \ins^T \ins -I)= \frac{k}{2}P\left(Q + \frac{k-2}{k}I\right)
\]
Since $P^2 = Q^2 = I$, then $P$ and $Q$ generate the dihedral group; that is to say, $\langle P, Q\rangle$ is a linear representation of the dihedral group. It is known that an indecomposable representation of this group over $\cx$ has dimension 1 or 2. Using this, we can compute the eigenvalues and multiplicities of elements of $\langle P, Q\rangle$, in particular, of $S^+(U)$. In \cite{Sze04}, Szegedy uses an observation of this flavour to find the spectrum of $U = PQ$. Here, we will use a similar decomposition of the Hilbert space and other linear algebra methods to explicitly determine the spectrum of $\splus$ in terms of the spectrum of the adjacency matrix. 

\begin{theorem}\label{thm:evals} If $G$ is a regular connected graph with valency $k \geq 2$ and $n$ vertices, then  $S^+(U(G))$ has eigenvalues as follows:
\begin{enumerate}[i)]
\item $k-1$ with multiplicity 1,
\item $\frac{\lambda \pm \sqrt{\lambda^2 - 4(k-1)}}{2}$ as $\lambda$ ranges over the eigenvalues of  $A$, the adjacency matrix of $G$, and $\lambda \neq k$,
\item 1 with multiplicity $\frac{n(k-2)}{2}+1$, and
\item $-1$ with multiplicity $\frac{n(k-2)}{2}$.
\end{enumerate} \end{theorem}

\proof For a matrix $M$, we write $\col(M)$ to denote the column space of $M$ and $\ker(M)$ to denote the kernel of $M$. Let $K = \col(\ins^T)+ \col(\outs^T)$ and let $L = \ker(\ins) \cap \ker(\outs)$. Observe that $K$ and $L$ are orthogonal complements of each other. Then $\re^{vk}$ is the direct sum of orthogonal subspaces $K$ and $L$. We will proceed by considering eigenvectors of 
$\splus$ in $K$ and in $L$ separately. 

For $K$, we will show that the eigenvectors of $\splus$ in $K$ lie in subspaces $C(\lambda)$ where $\lambda$ ranges over the eigenvalues of $A$. The eigenspace $C(k)$ has dimension $1$ while $C(\lambda)$ has dimension 2 for all $\lambda \neq k$. In $L$, we will show that all eigenvectors of $\splus$ have eigenvalue $\pm 1 $ and we will find the multiplicities of $\pm1$. 

First, we show that $K$ and $L$ are $\splus$-invariant. Since $L$ is the orthogonal complement of $K$, it suffices to check that $K$ is $\splus$-invariant. We obtain that:
\begin{equation}\label{SDins}
\splus \ins^T = k\outs^T - \outs^T = (k-1)\outs^T
\end{equation}
and 
\begin{equation}\label{SDouts}
\splus \outs^T  = \outs^TA - \ins^T .
\end{equation}
Hence, $K$ is $\splus$-invariant. 

We consider eigenvectors of $\splus$ in $K$. From equations \eqref{SDins} and \eqref{SDouts}, we obtain:
\begin{equation}\label{S2Douts} 
\splus^2 \outs^T = \splus (\outs^TA - \ins^T) 
= \splus \outs^TA - (k-1)\outs^T
\end{equation}

Let $\Zz$ be an eigenvector of $A$ with eigenvalue $\lambda$. Let $\Zy:= \outs^T \Zz$. Then, applying $\Zy$ to equation \eqref{S2Douts}, we obtain:
\[ \begin{split}
\splus^2 \Zy &= \splus^2 \outs^T \Zz \\
&= \splus \outs^TA\Zz - (k-1)\outs^T\Zz \\
&= \lambda \splus \Zy - (k-1)\Zy .
\end{split}
\]
Rearranging, we get
\begin{equation}\label{magic}
(\splus^2 - \lambda \splus + (k-1)I )\Zy = 0 .
\end{equation}

Let $C(\lambda) = \text{span}\{ \Zy, \splus \Zy \}$. By definition, $C(\lambda)$ has dimension at most 2 and is contained in $K$. For any vector $\Zx = \alpha\Zy + \beta \splus \Zy$ in $C(\lambda)$, we have that 
\[ \splus \Zx = \alpha \splus \Zy + \beta \splus^2 \Zy.
\] 
From equation \eqref{magic}, we can write $\splus^2 \Zy$ as a linear combination of $\splus \Zy $ and $\Zy$ and hence $\splus \Zx \in C(\lambda)$. Then, $C(\lambda)$ is $\splus$-invariant. If $C(\lambda)$ is 1-dimensional, then $\Zy$ is an eigenvector of $\splus$. Let $\theta$ be the corresponding eigenvalue. Then
\[ \begin{split}
\theta \Zy &= \splus \Zy \\
&= \splus \outs^T \Zz \\
&= ( \outs^TA - \ins^T) \Zz \\
&= \lambda \Zy - \ins^T \Zz 
\end{split}
\]
Then $(\theta - \lambda) \Zy = - \ins^T \Zz$ and $\Zz$ is in $\col(\ins^T) \cap \col(\outs^T)$. Then $\Zy$ is constant on arcs with a given head and on arcs with a given tail. Then $\Zy$ is constant on arcs of any component of $G$. Since $G$ is connected, $\Zy$ is the constant vector, which implies that $\Zz$ is a constant vector and $\lambda =k$. The eigenvalue of $\splus$ corresponding to $\Zy$ is $k-1$. 

Now suppose $C(\lambda)$ is 2-dimensional. Then, the minimum polynomial of $C(\lambda)$ is 
\[ t^2 - \lambda t + (k-1) = 0
\]
from \eqref{magic} and the eigenvalues are 
\[ \frac{\lambda \pm \sqrt{\lambda^2 - 4(k-1)} }{2} . \] 

These subspaces $C(\lambda)$ account for $2n -1$ eigenvalues of $\splus$. Since $\ins^T$ and $\outs^T$ are both $(nk) \times n$ matrices, $K$ has dimension at most $2n$. But, $\ins^T \mathbf{j} = \outs^T \mathbf{j} = \mathbf{j}$, where $\mathbf{j}$ is the all ones vector, since each row of both $\ins^T$ and $\outs^T$ has exactly one entry with value 1 and all other entries have value 0. Then, $K$ has dimension at most $2n -1$ and we have found all of the eigenvectors of $\splus$ in $K$. 

We will now find the remaining $n(k-2) + 1$ eigenvalues of $\splus$ over $L$. Let $\Zy$ be in $L$.  Then
\[ \begin{split}
\splus \Zy &= (\outs^T\ins - P) \Zy \\
&= \outs^T\ins\Zy  - P\Zy \\
&= - P\Zy .
\end{split}
\]
 If $\Zy$ is an eigenvector of $\splus$ with eigenvalue $\lambda$ and $\Zy$ is in $L$, then $\Zy$ is an eigenvector of $P$ with eigenvalue $-\lambda$. Since $P$ is a permutation matrix, $\lambda = \pm1$. 
 
To find the multiplicities we consider the sum of all the eigenvalues of $\splus$, which is equal to the trace of $\splus$. Observing that $P$ is a traceless matrix,
\[ \tr(\splus) = \tr(\outs^T\ins - P) = \tr(\outs^T\ins) = \tr(\ins\outs^T) = \tr(A) = 0 .
\]
The sum over all eigenvalues of $\splus$ should be 0. Let $\spec(A)$ be the set of eigenvalues of $A$. Consider the sum over the eigenvalues of eigenvectors of $K$:
\[ \begin{split}
&(k-1) + \sum_{\lambda \in \spec(A), \lambda \neq k}  \frac{\lambda \pm \sqrt{\lambda^2 - 4(k-1)} }{2}\\
&= (k-1) + \sum_{\lambda \in \spec(A), \lambda \neq k}  \lambda \\
&=-1  + \sum_{\lambda \in \spec(A)} \lambda \\
&= -1 .
\end{split}
\]
 
Then, the sum of the eigenvalue of the eigenvectors over $L$ is 1. So, 1 and $-1$ have multiplicities $\frac{n(k-2)}{2} + 1$ and $\frac{n(k-2)}{2}$, respectively. \qed

\section{Eigenvalues of $S^+(U^2)$}
\label{sec:su2}
We will show that $S^+(U^2) = (S^+(U))^2 + I$. Then, the eigenvalues of $S^+(U^2)$ are determined by the eigenvalues of $\splus$. The proof of the theorem will proceed by an analysis of which pairs of arcs give a negative entry in $U^2$.

\begin{theorem}\label{thm:su2} For any regular graph with valency $k$, if $k>2$ then $S^+(U^2) = S^+(U)^2 + I $.\end{theorem}
\proof Since $\outs^T \ins$ is the adjacency matrix of the line digraph of $G$, then $(\outs^T \ins)^2$ has the property that its $(j,i)$th entry counts the number of length two, directed walks in the line digraph of $G$. Observe that there is such a walk from $i$ to $j$ in $L(G)$ if and only if the head of $i$ is adjacent to the tail of $j$ in $G$. In particular, if there is a walk of length two from $i$ to $j$, there is only one such walk. Then, $(\outs^T \ins)^2$ is a $01$-matrix and is the support of $U^2$. We will find the required expression for $S^+(U^2)$ by subtracting from $(\outs^T \ins)^2$ the entries in $U^2$ which have negative value. 

We then proceed to look at the possible arrangements of $i$ and $j$ such that there is a length two, directed walk in $L(G)$ from $i$ to $j$, in Table \ref{3walks}.

\begin{table}[htdp]
\begin{center}
\begin{tabular}{|l|m{4.5cm}|m{3cm}|}
\hline
 & Directed walk of length 3 from $i$ to $j$ & Value of $(U^2)_{i,j}$\\
 \hline 
Case 1. & 
\begin{center}
\begin{tikzpicture}[decoration={
markings,
mark=at position .54 with {\arrow[black,thick]{>};}} ]
\draw[postaction={decorate}] (0,0) -- (1,0);
\draw[postaction={decorate}] (1,0) -- (2,0);
\draw[postaction={decorate}] (2,0) -- (3,0);
\filldraw[black]  (0,0) circle (1.5pt)
                           (1,0) circle (1.5pt)
                           (2,0) circle (1.5pt)
                           (3,0) circle (1.5pt);
\filldraw[white] (1.5, 0.3) circle (1pt);
\draw (0.5,0) node[anchor=north] {$i$}
           (2.5,0) node[anchor=north] {$j$};
\end{tikzpicture} 
\end{center}
 & \[\left( \frac{2}{k} \right)^2\] \\

\hline
Case 2. & 
\begin{center}
\begin{tikzpicture}[decoration={
markings,
mark=at position .54 with {\arrow[black,thick]{>};}} ]
\draw[postaction={decorate}] (0,0) -- (1,0);
\draw[postaction={decorate}] (1,0) -- (2,0);
\draw[postaction={decorate}]  (2,0) .. controls (1.3333,-0.4) and (0.6666, -0.4)..  (0,0);
\filldraw[black]  (0,0) circle (1.5pt)
                           (1,0) circle (1.5pt)
                           (2,0) circle (1.5pt);
\filldraw[white] (1.5, 0.3) circle (1pt);
\draw (0.5,0) node[anchor=south] {$i$}
           (1,-0.3) node[anchor=north] {$j$};
\end{tikzpicture} 
\end{center}
& \[\left( \frac{2}{k} \right)^2\] \\
\hline
Case 3. & 
\begin{center}
\begin{tikzpicture}[decoration={
markings,
mark=at position .54 with {\arrow[black,thick]{>};}} ]
\draw[postaction={decorate}] (0,0) -- (1,0);
\draw[postaction={decorate}] (1,0) .. controls (1.3333,0.3) and (1.6666, 0.3).. (2,0);
\draw[postaction={decorate}]  (2,0) .. controls (1.6666,-0.3) and (1.3333, -0.3)..  (1,0);
\filldraw[black]  (0,0) circle (1.5pt)
                           (1,0) circle (1.5pt)
                           (2,0) circle (1.5pt);
\filldraw[white] (1.5, 0.6) circle (1pt);
\draw (0.5,0) node[anchor=north] {$i$}
           (1.5,-0.2) node[anchor=north] {$j$};
\end{tikzpicture} 
\end{center}
&  \[\left( \frac{2}{k} \right)\left( \frac{2}{k} -1\right)\] \\
\hline
Case 4. & 
\begin{center}
\begin{tikzpicture}[decoration={
markings,
mark=at position .54 with {\arrow[black,thick]{>};}} ]
\draw[postaction={decorate}] (0,0) -- (1,0);
\draw[postaction={decorate}] (1,0) .. controls (1.3333,0.3) and (1.6666, 0.3).. (2,0);
\draw[postaction={decorate}]  (2,0) .. controls (1.6666,-0.3) and (1.3333, -0.3)..  (1,0);
\filldraw[black]  (0,0) circle (1.5pt)
                           (1,0) circle (1.5pt)
                           (2,0) circle (1.5pt);
\filldraw[white] (1.5, 0.6) circle (1pt);
\draw (1.5,0.3) node[anchor=south] {$i$}
           (0.5,0) node[anchor=north] {$j$};
\end{tikzpicture} 
\end{center}
&  \[\left( \frac{2}{k}-1 \right)\left( \frac{2}{k} \right)\]\\
\hline
Case 5. & 
\begin{center}
\begin{tikzpicture}[decoration={
markings,
mark=at position .54 with {\arrow[black,thick]{>};}} ]
\draw[postaction={decorate}] (1,0) -- (0,0);
\draw[postaction={decorate}] (0,0) .. controls (0.3333,0.5) and (0.6666, 0.5).. (1,0);
\draw[postaction={decorate}]  (0,0)  .. controls (0.3333,-0.5) and (0.6666, -0.5)..(1,0);
\filldraw[black]  (0,0) circle (1.5pt)
                           (1,0) circle (1.5pt);
%\filldraw[white] (1.5, 0.6) circle (1pt);
\draw (0.5,0.5) node[anchor=south] {$i$}
           (0.5,-0.5) node[anchor=north] {$j$};
\end{tikzpicture} 
\end{center}
&  \[\left( \frac{2}{k} -1 \right)^2\]\\
\hline
\end{tabular}
\end{center}
\caption{All possible pairs $i,j$ such that there is a length 2 walk in $L(G)$}
\label{3walks}
\end{table}%

We see that the only negative entries of $U^2$ occur for $i,j$ in Cases 3 and 4, when $k >2$. Then $(U^2)_{j.i}$ is negative when $i$ and $j$ share the same head but not the same tail and when $i$ and $j$ share the same tail but not the same head. Then,  
\[ \begin{split}
S^+(U^2) &= (\outs^T \ins)^2 -(\outs^T \outs -I) - (\ins^T \ins -I) \\
&= (\outs^T \ins)^2 -\outs^T \outs - \ins^T \ins + I +I \\
&= (\outs^T \ins)^2 -(\outs^T \ins)P - P(\outs^T \ins) + P^2 +I \\
&= (\outs^T \ins -P)^2 + I \\
&= S^+(U)^2 + I 
\end{split}
\]  \qed

The next theorem explicitly lists the eigenvalues of $S^+(U^2)$.

\begin{theorem} If $G$ is a regular connected graph with valency $k \geq 2$ and $n$ vertices, then  $S^+(U(G)^2)$ has eigenvalues as follows:
\begin{enumerate}[i)]
\item $k^2-2k + 2$ with multiplicity 1,
\item $\frac{\lambda^2 - 2k + 4}{2} \pm \frac{\lambda\sqrt{\lambda^2 - 4(k-1)}}{4}$ as $\lambda$ ranges over the eigenvalues of  $A$, the adjacency matrix of $G$, and $\lambda \neq k$ and
\item 2 with multiplicity $n(k-2)+1$.
\end{enumerate} \end{theorem}

\proof From Theorem \ref{thm:su2}, we get that $S^+(U^2) = (S^+(U))^2 + I$. Let $\Zy$ be an eigenvector of $\splus$ with eigenvalues $\theta$. Then,
$ S^+(U^2) \Zy = (\theta^2 + 1)\Zy $ and $\Zy$ is an eigenvector of $S^+(U^2)$ with eigenvalue $\theta^2 + 1$. The rest follow from the eigenvalues of $\splus$ found in Theorem \ref{thm:evals}. \qed

\section{Quantum Walk Algorithms for Graph Isomorphism}
\label{sec:QwalkGI}

The \textsl{Graph Isomorphism Problem} is the problem of deciding whether or not two given graphs are isomorphic. The algorithms of Shiau, Joynt and Coppersmith in \cite{SJC03},  Douglas and Wang in \cite{DW08}, and Gamble, Friesen, Zhou and Joynt in  \cite{GFZJ10} use the idea of evolving a quantum walk on a given pair of graphs and then comparing a permutation-invariant aspect of the states of the quantum walk on each graph. 

Both \cite{SJC03} and  \cite{GFZJ10} present algorithms based on a two-particle quantum walk.\footnote[1]{In the case of \cite{GFZJ10}, the particles are bosons.} Both procedures have been tested on large number of strongly regular graphs without finding a pair not distinguished by the procedure.  In \cite{S10}, Smith gives a family of graphs on which the procedure of Gamble et al.~\cite{GFZJ10} does not distinguish arbitrary graphs; in fact, he shows that $k$-boson quantum walks do not distinguish arbitrary graphs. However, the question of whether or not the procedure distinguishes all strongly regular graphs is still open. 

The quantum walk procedure of Douglas and Wang in \cite{DW08} has also been tested on classes of strongly regular graphs and of regular graphs, where all non-isomorphic graphs were distinguished. 

Finding a pair of non-isomorphic strongly regular graphs which are not distinguished by any of the three algorithms remains an open problem. Finding a pair of non-isomorphic strongly regular graphs which are not distinguished by the procedure of Emms et al.~is also an open problem. For work toward finding such a pair of graphs, see \cite{me10}. 

\bibliographystyle{plain}
%\bibliography{qwalk}

\end{document}